\input amstex
\input amsppt.sty

\magnification=1200
\voffset=-1.5cm

\NoBlackBoxes
\TagsOnRight
\NoRunningHeads

\def\p#1{\Cal P_#1}
\def\s#1{\frak S_#1}
\def\cs#1{\Bbb C[\frak S_#1]}
\def\a#1{\Cal A_#1}
\def\b#1{\Cal B_#1}
\def\bi{\b{\infty}}
\def\ai{\a{\infty}}
\def\deg{\operatorname{deg}}
\def\max{\operatorname{max}}

\def\res{\operatorname{Res}}

\def\t#1{\widetilde #1}
\def\ba#1{\overline #1}
\def\dgr{\downharpoonright}

\def\si{\sigma}
\def\ep{\epsilon}
\def\la{\lambda}
\def\al{\alpha}
\def\ga{\gamma}
\def\om{\omega}

\topmatter

\title
The algebra of conjugacy classes in symmetric groups,
and partial permutations
\endtitle
\author V.~Ivanov,\quad
$\boxed{\text{S.~Kerov}}$
\endauthor
\address
Moscow State University
\endaddress
\email
v\_n\_ivanov\@mtu-net.ru
\endemail
\address
Steklov Institute of Mathematics
at St.~Petersburg
\endaddress
\email
kerov\@pdmi.ras.ru.
\endemail

\thanks
{Partially supported by RFBR grants 98-01-00303
(the first author) and 99-01-00098 (the second author)}
\endthanks

\endtopmatter

\document

\subhead 1. Introduction \endsubhead
This article
was originally published in Russian
in
"Representation Theory, Dynamical Systems, Combinatorial and
Algorithmical Methods III" (A.~M.~Vershik, ed.), Zapiski Nauchnyh
Seminarov POMI {\bf 256} (1999), 95--120 (this text in Russian is
available via
{\tt http://www.pdmi.ras.ru/znsl/1999/v256.html}).

The main purpose of this note is to prove a convolution formula for conjugacy classes
in symmetric groups suggested in~\cite{7} (formula~(2.2), see also~\cite{8}).

Given a partition  $\rho\vdash r$ of a positive integer $r$, where $r\le
n$, denote by $\widetilde{\rho}=\rho\cup(1^{n-r})$
the partition of $n$ obtained from $\rho$ by adding an appropriate number of
unities.
Let $C_{\rho;n}$ be the sum of permutations $w\in\frak S_n$ of cycle type
$\widetilde{\rho}$. We also define ``normalized'' classes
$$
A_{\rho;n} = {n-r+m_1(\rho) \choose m_1(\rho)}\; C_{\rho;n},
\tag 1.1
$$
where $m_1(\rho)$ is the number of unities in the partition $\rho$.
We will show that there exist integer constants
$g_{\sigma,\tau}^\rho$ which define the convolution of normalized conjugacy
classes in the symmetric group
$\s n$:
$$
A_{\sigma;n}\,*\,A_{\tau;n} =
\sum_\rho g_{\sigma,\tau}^\rho\; A_{\rho;n}.
\tag 1.2
$$
Formula~(1.2) is valid for $n$ great enough. More exactly,
$n$ must be not less than $|\sigma|+|\tau|$, where $|\rho|$ denotes the
number such that $\rho$ is a partition of $|\rho|$; otherwise the summands in the
right-hand side which are not realized by permutations from
$\s n$ should be omitted. Note that the elements
$C_{\rho;n}$ do not form a basis in the algebra of central functions on
$\s n$: they may coincide for different partitions $\rho$.

Formula~(1.2) immediately implies an old result~\cite{5}:
$$
C_{\sigma;n}\,*\,C_{\tau;n} =
\sum_\rho q_{\sigma,\tau}^\rho(n)\; C_{\rho;n},
\tag 1.3
$$
where $\sigma,\,\tau,\,\rho$ run only over partitions without unity summands,
and the coefficients $q_{\sigma,\tau}^\rho(n)$ depend on $n$ in a polynomial
way.

To prove~(1.2), we introduce semigroups $\p n$ of partial permutations of
the set \linebreak $\{1,2,\ldots,n\}$. The semigroup algebras
$\Bbb C[\p n]$ are semi-simple and form a projective family with
limit $\b\infty=\varprojlim\Bbb C[\p n]$. The group of finite permutations
$\s\infty$ acts in $\b\infty$ by conjugations. Its orbits
$A_\rho$ are indexed by all partitions of positive integers, and form a basis
in the algebra of invariants $\ai=\b\infty^{\s\infty}$. The numbers
$g_{\sigma,\tau}^\rho$ arise as the multiplication structure constants
of the algebra $\ai$ in this basis.

We show that the algebra $\ai$ is naturally isomorphic to
the algebra of shifted symmetric functions $\Lambda^*$ introduced
in~\cite{3}. This isomorphism plays the same role for convolution of central
elements in the algebras $\Bbb C[\s n]$ as the characteristic mapping
$ch$ plays for the multiplication --- inducing of characters of
symmetric groups, see~\cite{9, I.7}. We also give examples of
filtrations on the algebra~$\ai$.

\subhead 2. The semigroup of partial permutations \endsubhead
Denote by $\Bbb P_n=\{1,\ldots,n\}$ a segment of positive integers, and
by $\frak S_n$ --- the group of all permutations of
$\Bbb P_n$. {\it A partial permutation} of the set $\Bbb P_n$ is a pair
$\alpha=(d,w)$ consisting of an arbitrary subset
$d\subset\Bbb P_n$ and an arbitrary bijection $w:d\to d$ of this subset.
The set $d$ will be referred to as {\it the support} of
$\alpha$. Denote by $\p n$ the set of all partial permutations of the set
$\Bbb P_n$.

Obviously, the number of elements in $\p n$ equals
$$
s_n = \sum_{k=0}^n {n \choose k}\, k! =
\sum_{k=0}^n (n\dgr k),
\tag 2.1
$$
where $(n\dgr k)=n(n-1)\ldots(n-k+1)$ is the falling factorial power.
These numbers satisfy a recurrence relation
$s_n= ns_{n-1} + 1$. Several first values are $s_0=1$,
$s_1=2$, $s_2=5$, $s_3=16$, $s_4=65$.

Given a partial permutation $(d,w)\in \p n$, denote by
$\widetilde{w}$ the permutation of the whole set $\Bbb
P_n$ coinciding with $w$ on $d$ and identical outside $d$.
The permutation $\widetilde{w}$ is well defined on all subsets containing
its support. This allows us to introduce a natural multiplication on the set
of partial permutations.

Given two partial permutations $(d_1,w_1)$, $(d_2,w_2)$, we define their
{\it product} as the pair $(d_1\cup
d_2,\,w_1w_2)$. With this multiplication, $\p n$ becomes a semigroup. The partial
permutation $(\varnothing,e_0)$, where $e_0$ is the trivial permutation of
the empty set $\varnothing$, is the unity in $\p n$, and it is a unique
invertible element in the semigroup
$\p n$.

\subhead 3. The semigroup algebra \endsubhead
Denote by $\b n=\Bbb C[\p n]$ the complex semigroup algebra of the semigroup
$\p n$. Let us check that this algebra is semi-simple and find its
irreducible representations.

Fix a subset $x\subset\Bbb P_n$ of size $|x|=k$ and denote by
$\frak S_x$ the group of permutations of the subset $x$.
It is easy to see that the formula
$$
\varphi_x(d,w) = \cases
\widetilde{w}, & \text{ if } d \subset x; \\
0 & \text{ otherwise}
\endcases
\tag 3.1
$$
defines a homomorphism of algebras $\b n\to\Bbb C[\s x]$.
This homomorphism is obviously surjective.

\proclaim{Lemma 3.1}
Let
$$
b = \sum_{k=0}^n\;\; \sum_{|d|=k}\; \sum_{w\in\s d}
b_{d,w}\, (d,w)
$$
be an element of the algebra $\b n$. Then the following conditions are
equivalent:
\roster
\item $\varphi_y(b)=0$ for all $y\subset x$;
\item $b_{d,w}=0$ for all $d\subset x$.
\endroster
\endproclaim

\demo{Proof}
Obviously, (2) implies~(1). To prove that (1) implies~(2), we use induction
on the size $k=|x|$ of the set $x$.

If $k=0$, we have $\varphi_\varnothing(b)=b_{\varnothing,e_0}\,e_0$, and
the statement is obvious.

Let $x\subset\Bbb P_n$ and let $v$ be a permutation of the set $x$.
Denote by $\ba d$ the set of non-fixed points of the permutation
$v$, and by $\ba v$ --- the restriction of $v$ on $\ba d$. Partial
permutations $(d,w)$ with $\varphi_x(d,w)=v$ are characterized by the
following conditions:
\roster
\item $\ba d \subset d \subset x$;
\item the restriction of $w$ on $\ba d$ coincides with $\ba v$;
\item all points of $d\setminus\ba d$ are fixed for $w$.
\endroster
Thus the coefficient of $v$ in the decomposition of $\varphi_x(b)$ equals
$$
\sum_{\ba d\subset d\subset x}\; \sum_{\ba w=\ba v}\; b_{d,w}\;,
$$
where $\ba w$ denotes the restriction of a permutation $w$ on
the set of its non-fixed points. The unique summand in this sum with
$d=x$ equals $b_{x,v}$. But by the induction hypothesis all summands with
$d\subsetneqq x$ are zero, hence this coefficient equals
$b_{x,v}$. The Lemma follows.
\qed\enddemo

\proclaim{Corollary 3.2}
The algebra $\b n$ is semi-simple, and it is isomorphic to the direct sum of the group
algebras of symmetric groups,
$$
\bigoplus_{x\subset\Bbb P_n} \varphi_x: \b n \cong
\bigoplus_{x\subset\Bbb P_n} \Bbb C[\s x].
\tag 3.2
$$
\endproclaim

\demo{Proof}
By Lemma~3.1, the homomorphism $\varphi=\bigoplus\varphi_x$ is injective,
and the dimension of the right-hand side equals
$$
\sum_{x\subset\Bbb P_n} |x|! =
\sum_{k=0}^n {n \choose k}\, k! = |\p n|
$$
and coincides with the dimension of $\b n$. Thus $\varphi$ is an isomorphism.
\qed\enddemo

Let
$$
\ep_d = \sum_{y:d\subset y\subset\Bbb P_n}
(-1)^{|y|-|d|}\; (y,e),
\tag 3.3
$$
where $e$ is the identity permutation. It is easy to see that
$\varphi_d(\varepsilon_d)=e\in\s d$ and
$\varphi_x(\varepsilon_d)=0$ for $x\ne d$. Thus the element
$\ep_d$ of the algebra $\b n$ is a central projection. Minimal central
projections are of the form
$\ep_d\delta_\lambda$, where $d\subset\Bbb P_n$ and $\delta_\lambda$
runs over minimal central projections in the algebra
$\Bbb C[\s d]$. The centre of the algebra $\b n$ is of the form
$$
Z(\b n) \cong \bigoplus_{d\subset\Bbb P_n}
Z(\Bbb C[\frak S_d]),
$$
where $Z(\Bbb C[\frak S_d])$ is the centre of the group algebra $\Bbb
C[\frak S_d]$. The dimension of the centre equals
$$
\dim Z(\b n) = \sum_{k=0}^n {n \choose k}\, p(k),
\tag 3.4
$$
where $p(k)$ is the number of partitions of $k$.

\subhead 4. Conjugacy classes in $\p n$ \endsubhead
The symmetric group $\frak S_n$ acts on the semigroup
$\p n$ by automorphisms $(d,w)\mapsto(vd,vwv^{-1})$. The orbits
of this action will be referred to as {\it conjugacy classes} in $\p n$.
It is obvious that two partial permutations are conjugate if and only if
the sizes of their supports coincide as well as their cycle types. Thus the
conjugacy classes $A_{\rho;n}\subset\p n$ are indexed by {\it partial
partitions} of $n$, i.e. by partitions $\rho\vdash r$ of any
integers $0\le r\le n$. In particular,
$A_{\varnothing;0}=\{(\varnothing,e)\}$.

Given a partial partition $\rho\vdash r\le n$, denote by
$\widetilde{\rho}=\rho\cup\{1^{n-r}\}$ the partition of $n$ obtained by
adding an appropriate number of unities. Let
$C_{\rho;n}$ be the conjugacy class in the group $\s n$ consisting of
permutations of cycle type $\widetilde{\rho}$. As usual, denote by
$m_k=m_k(\rho)$ the number of rows of length $k$ in the partition
$\rho$. The complement $\Bbb P_n\setminus d$ of the support of the partial
permutation $(d,w)$ contains $n-r$ points, thus the total number of
fixed points of $\widetilde{w}$ equals $n-r+m_1(\rho)$.

Denote by $\psi:(d,w)\mapsto\widetilde{w}$ the homomorphism of forgetting
the support of a partial permutation.

Let a permutation $v\in\s n$ have cycle type
$\widetilde{\rho}$. The set $\psi^{-1}(v)\cap A_{\rho;n}$
consists exactly of partial permutations $(d,w)$ such that the support
$d$ contains all non-fixed points of the permutation $v$. In the set
$d$, one may arbitrarily choose $m_1(\rho)$ fixed points from the total
number of fixed points of the permutation $v$ which is equal to
$n-r+m_1(\rho)$. Hence the numbers of elements in the class
$A_{\rho;n}$ and in the conjugacy class $C_{\rho;n}$ of the group $\s n$
are related by the formula
$$
|A_{\rho;n}| =
{n-r+m_1(\rho) \choose m_1(\rho)}\, |C_{\rho;n}|.
\tag 4.1
$$

The action of the symmetric group $\s n$ on $\p n$ can be continued by
linearity to an action of $\s n$ on the algebra $\b n$. Denote by
$\a n=\b n^{\s n}$ the subalgebra of invariant elements for this action.

The homomorphism $\psi:\p n\to\s n$ can also be continued to a surjective
homomorphism of algebras $\psi:\b n\to\cs n$. It commutes with the action of
the group $\s n$ by conjugations on the algebras
$\b n$ and $\cs n$. Thus $\psi(\a n)=Z(\cs n)$, where $Z(\cs n)$ is the
centre of the group algebra $\cs n$.

Let us identify the conjugacy class $A_{\rho;n}$ with the element
$$
A_{\rho;n}=\sum_{(d,w)\in A_{\rho;n}}(d,w)
\tag 4.2
$$
of the algebra $\b n$. In particular, if $|\rho|>n$, then $A_{\rho;n}=0$.
It follows from our definitions that the elements $A_{\rho;n}$, where
$|\rho|\le n$, form a linear basis of the algebra $\a n$. It is clear that
$$
\psi(A_{\rho;n}) =
{n-r+m_1(\rho) \choose m_1(\rho)}\, C_{\rho;n}.
\tag 4.3
$$

In Sect.~12 we construct all irreducible representations
$\pi_{x,\la}$ of the algebra $\b n$. Note that in any irreducible
representation $\pi_{x,\la}$ of the algebra $\b n$ the element
$A_{\rho;n}$ acts as a scalar operator. Thus
$A_{\rho;n}\in Z(\b n)$, and the algebra $\a n$ lies in the centre
$Z(\b n)$. This inclusion is strict for $n\ge2$.

\subhead 5. Algebras $\bi$ and $\ai$ \endsubhead
Let $m\le n$. We introduce a mapping $\theta_m:\b n\to\b m$
by the formula
$$
\theta_m(d,w) = \cases
(d,w), & \text{if } d\subset \Bbb P_m,\\
0 & \text{otherwise}.
\endcases
\tag 5.1
$$
The mapping $\theta_m$ is a homomorphism of algebras and it commutes with the
action of the group $\s m$ on $\b n$ and on $\b m$.
Hence $\theta_m(\a n)=\a m$.

Define {\it the degree} of a partial permutation $(d,w)\in\p n$
as $\deg(d,w)=|d|$. Given $b=\sum_{\al\in\p
n}b_\alpha\,\alpha\in\b n$, let $\deg b=\max\deg\alpha$, where the maximum
is over all $\alpha$ with $b_\alpha\ne0$. The function
$\deg$ defines a filtration on the algebra $\b n$. Note that
$\deg(\theta_m(b))\le\deg b$ for all $b\in\b n$.

Denote by $\bi$ the projective limit of the algebras $\b n$ with respect to
the morphisms $\theta_n$, and by $\ai$ --- the projective limit of the
algebras $\a n$. Both limits are taken in the category of filtered
algebras.

Let $\s{\infty}$ be the infinite symmetric group, i.e. the group of finite
permutations of positive integers. The group
$\s{\infty}$ acts naturally on  $\bi$, and
$\ai=\bi^{\s{\infty}}$ is the subalgebra of invariants for this action.

\subhead 6. Structure constants of the algebra $\ai$ \endsubhead
Denote by $\theta_n$ the natural homomorphism
$\theta_n:\bi\to\b n$ as well as its restriction on $\ai$.
The natural inclusion of algebras $i_n:\b n\to\bi$ accords with the
projection $\theta_n$: \quad $\theta_n\circ i_n=id_{\b n}$.
It is convenient to write elements of $\bi$ as formal infinite sums,
$$
b = \sum_{n=0}^\infty\;\; \sum_{|d|=n}\; \sum_{w\in\s d}
b_{d,w}\, (d,w).
\tag 6.1
$$

Given a partition $\rho\vdash r$, let $A_\rho=\sum
(d,w)$; the sum extends to partial permutations
$(d,w)\in\p \infty$ such that $|d|=r$ and $w$ has cycle type $\rho$. The
elements $A_\rho$, where $\rho$ runs over all partitions, form a linear basis
in $\ai$.

Denote by $g_{\si,\tau}^\rho$ the structure constants of the algebra
$\ai$ in the basis $\{A_\rho\}$,
$$
A_\si\; A_\tau =
\sum_\rho g_{\si,\tau}^\rho\, A_\rho.
\tag 6.2
$$
Note that $\theta_n(A_\rho)=A_{\rho;n}$, where $A_{\rho;n}$ is the element
of the algebra $\a n$ introduced in Sect.~4. Since
$\theta_n:\ai\to\a n$ is a homomorphism, we obtain the following statement.

\proclaim{Proposition 6.1}
$$
A_{\si;n}\; A_{\tau;n} =
\sum_{\rho}g_{\si,\tau}^\rho\, A_{\rho;n}.
\tag 6.3
$$
\endproclaim

For $|\rho|\le n$, by definition $A_{\rho;n}=0$. Let us illustrate
Proposition~6.1 by an example. The simplest non-trivial multiplication
formula in the algebra $\ai$ is
$$
A_{(2)}\,A_{(2)}=A_{(1^2)}+3A_{(3)}+2A_{(2^2)}.
$$
In the algebras $\a 2$, $\a 3$, $\a 4$ we have
$$
\align
A_{(2);2}\;A_{(2);2} &= A_{(1^2);2}\\
A_{(2);3}\;A_{(2);3} &= A_{(1^2);3}+3A_{(3);3}\\
A_{(2);4}\;A_{(2);4} &= A_{(1^2);4}+3A_{(3);4}+2A_{(2^2);4}.
\endalign
$$

Let us give a useful combinatorial interpretation of the structure constants
$g_{\si,\tau}^\rho$.

\proclaim{Proposition 6.2}
Given a partition $\rho\vdash r$, let $d_\rho=\Bbb P_r$ and let
$$
w_\rho=(1,\dots,\rho_1)(\rho_1+1,\dots,\rho_1+\rho_2)\dots
(|\rho|-\rho_{\ell(\rho)}+1,\dots,|\rho|)
$$
be a fixed permutation of the set $d_\rho$ of cycle type
$\rho$. Consider the set $G_{\si,\tau}^\rho(n)$ of pairs
$\big((d_1,w_1),(d_2,w_2)\big)\in\p n\times\p n$ such that
\roster
\item\quad $(d_1,w_1)\in A_{\si;n},\quad (d_2,w_2)\in A_{\tau;n}$;
\item\quad $d_1\cup d_2=d_\rho$,\quad $w_1\,w_2=w_\rho$.
\endroster
Then for $n\ge r$ the number of elements $|G_{\si,\tau}^\rho(n)|$
equals $g_{\si,\tau}^\rho$.
\endproclaim

\demo{Proof}
The partial permutation $(d_\rho,\,w_\rho)$ belongs to the class
$A_\rho$. By definition, the number $|G_{\si,\tau}^\rho(n)|$ is the
coefficient of the element $(d_\rho,\,w_{\rho})$ in the product
$A_{\si;n}\,A_{\tau;n}$. If $n\ge |\rho|$, then $A_{\rho;n}\ne
0$ and using Proposition~6.1 we obtain
$|G_{\si,\tau}^\rho(n)|=g_{\si,\tau}^\rho$.
\qed\enddemo

\proclaim{Proposition 6.3}
If $g_{\si,\tau}^\rho\ne0$, then $|\rho|\le |\si|+|\tau|$.
\endproclaim

\demo{Proof}
If $g_{\si,\tau}^\rho\ne 0$, then it follows from Proposition~6.2 that
$G_{\si,\tau}^\rho(|\rho|)\ne\varnothing$. Hence there exist sets
$d_1,d_2$ with $|d_1|=|\si|,\;
|d_2|=|\tau|,\; |d_1\cup d_2|=|\rho|$. Thus $|\rho|\le
|\si|+|\tau|$.
\qed\enddemo

\remark{Remark 6.4}
Denote by $\si\cup\tau$ the partition whose parts are obtained by uniting
all parts of partitions $\si$ and $\tau$. It follows from Proposition~6.2
that $\rho=\si\cup\tau$ is the unique partition with
$|\rho|=|\si|+|\tau|$ and $g_{\si,\tau}^\rho\ne0$, the coefficient being
equal to
$$
g_{\si,\;\tau}^{\si\cup\tau}=\prod_{k\ge 1}
{m_k(\si)+m_k(\tau)\choose m_k(\si)}.
$$
\endremark

\subhead 7. Convolutions of conjugacy classes in symmetric groups\endsubhead
Given an arbitrary partition $\rho\vdash r$, denote by
$\ba\rho$ the partition obtained from
 $\rho$ by removing all its unity parts (if they existed). Thus
$m_1(\ba\rho)=0$ and $\rho=\ba\rho\cup 1^{r-|\ba{\rho}|}$.

Recall that we have selected
in the centre of the group algebra of the symmetric group
$Z(\cs n)$ the elements $C_{\rho;n}$ indexed by partitions
$\rho\vdash r$ of any numbers $0\le
r\le n$,
$$
C_{\rho;n} = \sum_{w{\text{
\ is of type\ }}\t{\rho}} w.
\tag 7.1
$$
In particular, $C_{\rho;n}=0$, if $|\rho|>n$. Assume now that
$|\si|,|\tau|\le n$; then $C_{\si;n}=C_{\tau;n}$ if and only if
$\overline\si=\overline\tau$. Let us say that a partition
$\rho\vdash r$ is {\it proper}, if $m_1(\rho)=0$. Denote the set of all
proper partitions by $\overline{\Bbb Y}$. The elements
$\{C_{\rho;n}:\;\rho\in\overline{\Bbb Y},\;|\rho|\le n\}$
form a linear basis in $Z(\cs n)$.

\proclaim{Theorem 7.1}
Given a partition $\rho\vdash r\le n$, consider the images
$\psi(A_{\rho;n})$
$$
\psi(A_{\rho;n}) = {n-r+m_1(\rho)\choose m_1(\rho)}\;C_{\rho;n}
\tag 7.2
$$
of the elements $A_{\rho;n}\in\a n$ in the centre of the group algebra of the
symmetric group $\frak S_n$ under the ``forgetting support'' mapping from
Sect.~4. Then

{\rm a)} for every $n$, the following equality holds,
$$
\psi\big(A_{\si;n}\big)\; \psi\big(A_{\tau;n}\big) =
\sum_{\rho} g_{\si,\tau}^\rho\, \psi\big(A_{\rho;n}\big),
\tag 7.3
$$
where $g_{\si,\tau}^\rho$ are the structure constants of the algebra
$\ai$ which do not depend on $n$.

{\rm b)} for $n\ge|\si|+|\tau|$, the sum in the right-hand side of the
equality~$(7.3)$ is stable, i.e. the non-zero summands are indexed by the same
partitions $\rho$.
\endproclaim

\demo{Proof}
a) Equality~(7.3) follows from the fact that $\psi$ is a homomorphism and
from Proposition~6.1.

b) The collection of non-zero summands in the right-hand side of equality~(7.3)
is indexed by partitions  $\rho$ such that $|\rho|\le n$ and
$g_{\si,\tau}^\rho\ne 0$. By Proposition~6.3, this collection is fixed for $n\ge
|\si|+|\tau|$.
\qed\enddemo

\proclaim{Proposition 7.2}
Let numbers $h_{\si,\tau}^\rho$ satisfy
$$
\psi\big(A_{\si;n}\big)\; \psi\big(A_{\tau;n}\big) =
\sum_{\rho} h_{\si,\tau}^\rho\, \psi\big(A_{\rho;n}\big)
$$
for every $n$. Then $h_{\si,\tau}^\rho=g_{\si,\tau}^\rho$ for all
$\rho,\si,\tau$.
\endproclaim

\demo{Proof} Suppose the contrary.
Choose arbitrary partitions
$\si$ and $\tau$. Let $\rho$ be a partition such that
$h_{\si,\tau}^\rho\ne g_{\si,\tau}^\rho$ and
$h_{\si,\tau}^\nu=g_{\si,\tau}^\nu$ for any partitions $\nu$
with $|\nu|<|\rho|$. Let $|\rho|=n$; then
$$
0 = \sum_\nu (g_{\si,\tau}^\nu-h_{\si,\tau}^\nu)\,\psi(A_{\nu;n}) =
\sum_{|\nu|=n} (g_{\si,\tau}^\nu-h_{\si,\tau}^\nu)\,\psi(A_{\nu;n}).
$$
On the other hand, the set $\{\psi(A_{\nu;n}):\,\,|\nu|=n\}$
is linearly independent. The obtained contradiction proves the Proposition.
\qed\enddemo

Thus $\{\psi(A_{\rho;n})\}$ is a family of elements proportional to
conjugacy classes (with proportionality coefficients depending on
$\rho$ and on $n$), and the multiplication structure constants do not depend on
$n$. A similar family was introduced earlier in~\cite{7, 8}.
The notion of such elements is also close to~\cite{1}.

Let us come back to multiplication of conjugacy classes $C_{\rho;n}$.

\proclaim{Proposition 7.3}
Let $\sigma,\tau,\rho\in\overline{\Bbb Y}$ be proper partitions
(i.e. without unity parts).
We define polynomials $q_{\si,\tau}^\rho(n)$ as
$$
q_{\si,\tau}^\rho(n) = \sum_{k\ge0}
g_{\si,\;\tau}^{\rho\cup (1^k)}{n-|\rho|\choose k}.
\tag 7.4
$$
Then
$$
C_{\si;n} * C_{\tau;n} =
\sum_{|\rho|\le n,\, m_1(\rho)=0}
q_{\si,\tau}^\rho(n)\, C_{\rho;n}
\tag 7.5
$$
for all $\sigma,\tau\in\overline{\Bbb Y}$.
\endproclaim

The fact that the coefficients $q_{\si,\tau}^\rho(n)$ in~(7.5)
are polynomials on $n$ assuming only integer values in integer points was
first obtained in~\cite{5, Theorem~2.2}.

\demo{Proof}
For a proper partition $\sigma$, formula~(4.3) becomes simpler and reduces to
$\psi(A_{\sigma;n})=C_{\sigma;n}$. Hence
$$
\align
C_{\si;n}\; C_{\tau;n} =
\psi(A_{\si;n})\; \psi(A_{\tau;n}) =
&\sum_{|\rho|\le n} g_{\si,\tau}^\rho\,
\psi(A_{\rho;n}) =\\
=&\sum_{|\rho|\le n} g_{\si,\tau}^\rho
{n-|\rho|+m_1(\rho) \choose m_1(\rho)}
C_{\rho;n}.
\endalign
$$
Collecting similar summands of the form
$C_{\rho;n}=C_{\ba\rho\cup(1^k);n}$ for $k\ge0$, we obtain formula~(7.4)
for coefficients in~(7.5).
\qed\enddemo

It follows from Remark~6.4 that $\deg
q_{\si,\;\tau}^{\si\cup\tau}=0$.

\remark{Remark 7.4}
The sum participating in formula~(7.4) may contain more than one non-zero
summand. Consider, for example, $\si=\tau=(3)$. Formulae from Sect.~11 show
that both
$g_{\si,\tau}^{(3)}\ne 0$
and
$g_{\si,\tau}^{(3,1)}\ne 0$.
\endremark

\remark{Remark 7.5}
Let $\Phi$ be the ring of polynomials in one variable assuming only integer
values in integer points, and let $\Phi_{\Bbb C}$ be the ring of polynomials of
the form $\t P(t)=cP(t)$, where $P\in\Phi$ and $c\in\Bbb C$. In~\cite{5}
the authors deal with a $\Phi$-algebra $\Cal K$ freely generated over
$\Phi$ by a basis $C_\rho$, where
$\rho\in\overline{\Bbb Y}$ runs over proper partitions (without unity parts).
A formal relation of our algebra $\ai$ and the algebra
$\Cal K$ is given by the fact that the mapping
$$
A_\rho\mapsto {n-|\rho|+m_1(\rho) \choose m_1(\rho)} C_{\ba\rho}
$$
defines an epimorphism of $\Phi_{\Bbb C}$-algebras
$\ai\bigotimes_{\Bbb C}\Phi_{\Bbb C}$ and $\Cal
K\bigotimes_\Phi\Phi_{\Bbb C}$.
\endremark

\subhead 8. The semigroup of fillings \endsubhead
Given a partition $\rho\vdash r$, denote by
$$
z_\rho=\prod_{i\ge 1}i^{m_i(\rho)}m_i(\rho)!
$$
the size of the centralizer of a permutation of cycle type $\rho$.
Note that
$$
{z_\rho \over m_1(\rho)!} = z_{\ba\rho} =
{z_{\t\rho} \over (n-r+m_1(\rho))!}\;.
\tag 8.1
$$

Along with the introduced above elements $A_\rho$,
we consider a related basis in the algebra $\ai$ consisting of the elements
$a_\rho=z_\rho\,A_\rho$. Let $f_{\sigma,\tau}^\rho$ be the multiplication
structure constants in this new basis,
$$
a_{\si;n}\; a_{\tau;n} =
\sum_{\rho} f_{\si,\tau}^\rho\, a_{\rho;n}.
\tag 8.2
$$
They are related to the constants $g_{\sigma,\tau}^\rho$ by an obvious
formula
$$
f_{\si,\tau}^\rho = \frac{z_\si z_\tau}{z_\rho}\;g_{\si,\tau}^\rho.
\tag 8.3
$$
We will show that all numbers $f_{\sigma,\tau}^\rho$ are non-negative
integers and give their combinatorial interpretation similar to
Proposition~6.2.

\smallskip\noindent
{\bf Definition.}
Let $\lambda$ be a Young diagram with $k$ boxes, and $d$ be the set of
$k$ distinct positive integers. Any bijection
$R:\lambda\to d$ will be referred to as a {\it filling} of shape
$\lambda$ and weight $d$, and the set $d$ will be called {\it the support}
of the filling $R$. Define a permutation $w_R\in\s d$ by declaring
the rows of the filling $R$ to be cycles, so that the cycle type of $w_R$ is
$\lambda$. Obviously, the number of fillings of shape $\lambda$ and a
fixed weight equals $n!=z_\lambda\,|C_\lambda|$.

\smallskip
For example, given a filling
$$
R = \pmatrix 4&3&1\\9&2&7\\6&5 \endpmatrix
$$
of the diagram $\lambda=(3,3,2)$ with support $d=\{1,2,3,4,5,6,7,9\}$, we
have $w_R=(1,4,3)$ $(2,7,9)$ $(5,6)$.

We define {\it the convolution $R=S*T$ of fillings} $S,T$ by the following
rules.

Let $d_R=d_S\cup d_T$ be the union of the supports of fillings
$S,T$. Order the set $d_R$ by reading first the elements of $S$ and then the
elements of $T$ from left to  right along each row and from top to bottom.
Repeating
elements of the filling $T$ are to be ignored.

The first element of $d_R$ is the first element $s$ of the first row of the
filling $S$. Form a row of the filling $R$, which we want to define, as the cycle of the
product of permutations $\t w_S$ and $\t w_T$ containing
$s$ and beginning with $s$. The tilde means that the
permutations $w_S$, $w_T$ are trivially continued to be
defined on the set
$d_R$: all points outside the former domain
are assumed fixed.

Passing to constructing other rows of the filling $R$, assume that a part of
rows of $R$ is already constructed. If the elements of the set
$d_R$ are not all used, denote by $s$ the first of the remained elements
(in the above-mentioned order). The next row of $R$ is the cycle of the product of
permutations $\t w_S$ and $\t w_T$ containing
$s$. The first element of the row is chosen to be $s$.

The lengths of rows formed according to these rules do not necessarily
decrease. Let us reorder the rows by decreasing of lengths without changing
the respective order of rows of equal length. The obtained filling $R$ is
the convolution of fillings $S$ and $T$.

We illustrate the definition of the convolution by an example. Let
$$
S = \pmatrix 3&4&5&6&9\\2&1&7 \endpmatrix; \qquad
T = \pmatrix 4&3&2\\1&9&6\\8 \endpmatrix.
$$
Then $d_R=\Bbb P_9$, and the product of the permutations $\t w_S$ and $\t
w_T$ is equal to \linebreak $(1,3)(2,5,6,7)(4)(8)(9)$. Hence
$$
S*T = \pmatrix 3&1\\4\\5&6&7&2\\9\\8 \endpmatrix =
\pmatrix 5&6&7&2\\3&1\\4\\9\\8 \endpmatrix.
$$

\proclaim{Proposition 8.1}
Fix Young diagrams $\sigma$, $\tau$. Given a partition
$\rho\vdash r$, put $d_\rho=\Bbb P_r$ and let $R_\rho$ be the canonical
filling of the Young diagram $\rho$ in which the boxes are indexed with the
numbers $1,2,\ldots,r$ successively from left to right and
from top to bottom.
Consider the set $F_{\si,\tau}^\rho$ of pairs of fillings $(S,T)$ such that
\roster
\item\quad $S$ is of shape $\sigma$, and $T$ is of shape $\tau$;
\item\quad $S*T=R_\rho$.
\endroster
Then the number of elements $|F_{\si,\tau}^\rho|$ equals
$f_{\si,\tau}^\rho$. In particular, the structure constants in~$(8.2)$
are non-negative integers.
\endproclaim

\demo{Proof}
The set $\Cal F_r$ of fillings $R$ with support $d_R\subset\Bbb
P_r$ forms a semigroup with respect to the introduced above convolution
operation. Associating a partial permutation
$(d_R,w_R)$ with a filling $R$, we obtain an epimorphism of semigroups
$u:\Cal F_r\to\Cal
P_r$.

The group $\s r$ acts in an obvious way by automorphisms of the semigroup
$\Cal F_r$, and the homomorphism $u$ is equivariant under this action. Each
partial permutation $(d,w)\in\Cal P_r$ of cycle type
$\rho$ has exactly $z_\rho$ inverse images in $\Cal
F_r$. Thus the multiplication structure constants
$|F_{\si,\tau}^\rho|$ and $g_{\si,\tau}^\rho$ are related by a formula of
type~(8.3), and the Proposition follows.
\qed\enddemo

One can easily deduce from Proposition~8.1 formulae for coefficients
$f_{\sigma,\tau}^\rho$ in the simplest cases.

\proclaim{Corollary 8.2}
If $\rho=\sigma\cup\tau$, then $f_{\sigma,\tau}^\rho=1$.
\endproclaim

\proclaim{Corollary 8.3}
Let a Young diagram $\rho$ be obtained from Young diagrams $\sigma$,
$\tau$ by $(1)$ replacing a row of $\sigma$ of length $i$ and a row of $\tau$
of length $j$ by a row of length  $i+j-1$, and $(2)$ uniting the remained rows
(with subsequent ordering by decreasing). Let
$m_i(\sigma)$ be the multiplicity of rows of length $i$ in $\sigma$, and
$m_j(\tau)$ be the multiplicity of rows of length $j$ in $\tau$. Then
$f_{\sigma,\tau}^\rho=im_i(\sigma)\,jm_j(\tau)$.
\endproclaim

\subhead 9. Isomorphism of the algebra $\ai$ and the algebra of shifted
symmetric functions\endsubhead
The algebra $\Lambda^*$ of shifted symmetric functions was introduced and
studied in~\cite{10, 8, 3}. In this section we establish an
isomorphism of this algebra with the algebra $\ai$ and indicate the
elements of $\ai$ corresponding to the shifted Schur functions
$s^*_\lambda$ and shifted analogues of the Newton power sums
$p^\#_\rho$ introduced in~\cite{3,~(1.6) and~(14.9)}.

The algebra $\Lambda^*$ is defined as follows. Denote
by $\Lambda^*(n)$ the algebra of polynomials with complex coefficients in
$x_1,\,\dots,x_n$ that become symmetric in new variables
$x_i'=x_i-i$. The algebra
$\Lambda^*(n)$ is filtered by the degree of polynomials. The specification
$x_{n+1}=0$ defines a homomorphism of filtered algebras
$\Lambda^*(n+1)\to\Lambda^*(n)$. Denote by $\Lambda^*$ the projective limit
of the algebras $\Lambda^*(n)$ with respect to these homomorphisms
(in  the category of filtered algebras).
The  algebra $\Lambda^*$ is called {\it the algebra of shifted symmetric
functions.} The ring $\Lambda^*$ can also be defined over $\Bbb
Z$, but for our purposes it is more convenient to assume
that $\Lambda^*$ is an algebra over the field $\Bbb C$.

Given an element $f\in\Lambda^*$ and a partition $\la$, we denote by
$f(\la)$ the value $f(\la_1,\,\dots,\la_{\ell(\la)})$.
Elements of the algebra $\Lambda^*$ are uniquely defined by their values on
partitions.

The key point of the paper~\cite{3} is the basis of {\it the shifted
Schur functions} $\{s_\mu^*\}$ of the space $\Lambda^*$ indexed by
partitions $\mu\in\Bbb Y$. The paper~\cite{3} contains explicit formulae and
many other remarkable facts for these functions, but we use below only the
following two properties of the functions
$s_\mu^*$. If $\rho\vdash r\le n$, then
$$
\sum_{\mu\vdash r}\frac{s_\mu^*(\la)}{(n\dgr r)}\,\chi_\rho^\mu=
\frac{1}{\dim\lambda} \sum_{\mu\vdash r}
\dim(\lambda/\mu)\, \chi_\rho^\mu =
\frac{1}{\dim\lambda} \chi^\lambda_{\tilde\rho}
\tag 9.1
$$
according to~\cite{3, Theorem~7.1}. And if $\mu\vdash r>|\lambda|$,
then
$$
s_\mu^*(\la)=0
\tag 9.2
$$
by~\cite{3, Theorem~3.1}. Given a partition $\lambda\vdash n$, we denote by
$\chi^\la$ the irreducible character of the symmetric group
$\frak S_n$. If $\rho$ is a partition of a number $r\le n$, then
$\chi_{\t\rho}^\la$ is the value of the character $\chi^\la$ on an element
of cycle type $\t\rho$.

Following~\cite{3, \S14.2} (see also~\cite{8}), we introduce another basis
$\{p_\rho^\#\}$ in $\Lambda^*$ by the formula
$$
p_\rho^\# = \sum_{\mu\vdash r} \chi^\mu_\rho\, s^*_\mu.
$$
Note that $\deg p_\rho^\#=|\rho|$. It follows from~(9.1) and~(9.2) that
$$
{p_\rho^\#(\la) \over (n\dgr r)} =
\cases \frac{1}{\dim\la}\; \chi_{\t\rho}^\la,
&\text{ if } n\ge r,\\
0, &\text{ otherwise}.
\endcases
\tag 9.3
$$

\proclaim{Theorem 9.1}
The linear mapping $F:\ai\to\Lambda^*$ defined on the basis elements of
$A_\rho$ by the formula
$$
F(A_\rho) = {p_\rho^\# \over z_\rho}
\tag 9.4
$$
is an isomorphism of algebras $\ai$ and $\Lambda^*$.
\endproclaim

\demo{Proof}
Let $\rho$ be a partition of a number $r$ and $m_1=m_1(\rho)$. According to~(4.3),
the image of the element $A_\rho\in\ai$ in the centre of the group algebra
$Z(\Bbb C[\frak S_n])$ equals
$$
(\psi\circ \theta_n)(A_\rho) = {n-r+m_1 \choose m_1} C_{\rho;n}.
$$
If $r\le n$, then
$$
\chi^\lambda\big(\psi\circ \theta_n(A_\rho)\big) =
{n-r+m_1 \choose m_1}\;{n! \over z_{\t\rho}}\;\chi^\lambda_{\t\rho}
= {(n \dgr r) \over z_\rho}\;\chi^\lambda_{\t\rho},
$$
thus formula~(9.3) implies
$$
{1 \over z_\rho} p_\rho^\#(\lambda) =
{(n \dgr r) \over z_\rho}\;{1 \over \dim\lambda}\;
\chi^\lambda_{\t\rho} =
{1 \over \dim\lambda}\;\chi^\lambda\big(\psi\circ
\theta_n(a_\rho)\big).
$$
If $r>n$, the last formula is valid too, since both sides
are zero in view of~(9.2). Thus the mapping
$F:\ai\to\Lambda^*$ defined in~(9.4) may be defined
by an equivalent formula
$$
\big(F(a)\big)(\la) =
{1 \over \dim\la} \chi^\la (\psi\circ \theta_{|\la|}(a)); \qquad
a \in \ai.
\tag 9.5
$$
The mapping $\psi$ is a homomorphism of the algebra $\a n$ onto $Z(\Bbb
C[\frak S_n])$, and irreducible normalized characters
$\chi^\lambda/\dim\lambda$ of the group $\frak S_n$ define homomorphisms
$Z(\Bbb C[\frak S_n])\to\Bbb C$. Hence,
$$
F(ab)(\la)=(F(a)\, F(b))(\la)
\tag 9.6
$$
for all $a,b\in\ai$, and formulae~(9.4),~(9.5) define an isomorphism of the
algebra $\ai$ onto the algebra $\Lambda^*$. The Theorem follows.
\qed\enddemo

\proclaim{Proposition 9.2}
Associate with a partition $\mu\vdash m$ an element
$$
x_\mu=\sum_{|d|=m}\; \sum_{w\in\s d}
\chi^\mu(w)\, (d,w)
\tag 9.7
$$
of the algebra $\ai$. Then $F(x_\mu)=s_\mu^*$.
\endproclaim

\demo{Proof}
If $|\mu|>|\lambda|$, then $F(x_\mu)(\la)=0=s_\mu^*(\lambda)$.
If $|\mu|\le|\lambda|$, then
$$
\psi\circ \theta_n(x_\mu) =
\sum\Sb |d|=m \\ d\subset\Bbb P_n\endSb\;
\sum_{w\in\s d} \psi(d,w) =
\sum\Sb |d|=m \\ d\subset\Bbb P_n\endSb\;
\sum_{w\in\s d} \t w.
$$
{}From~(9.1) and the
orthogonality relations for irreducible characters of the
symmetric group $\s m$, we obtain a chain of equalities
$$
\aligned
F(x_\mu)(\lambda) &=
\sum\Sb |d|=m \\ d\subset\Bbb P_n\endSb
\sum_{w\in\s d} \chi^\mu(w)\;
\frac{\chi^\lambda(\t w)}{\dim\lambda} = \\ &=
{n \choose m}\, \sum_{w\in\frak S_m}
\chi^\mu(w) \; \sum_{\nu\vdash m}
\chi^\nu(w) \frac{s_\nu^*(\la)}{(n\dgr m)} = \\ &=
{n \choose m}\, \sum_{\nu\vdash m} \frac{s_\nu^*(\la)}{(n\dgr m)}
\sum_{w\in\frak S_m} \chi^\mu(w) \, \chi^\nu(w) = \\ &=
{n \choose m}\, \frac{m!\;s_\mu^*(\la)}{(n \dgr m)} = s_\mu^*(\la).
\endaligned
$$
The Proposition follows.
\qed\enddemo

\subhead 10. Filtrations of the algebra $\ai$ \endsubhead
There is an obvious filtration on the algebra $\ai$,
$$
\deg_1(A_\rho)=|\rho|.
\tag 10.1
$$
In the decomposition of the convolution $A_\si\,A_\tau$, there is the unique
summand $A_{\si\cup\tau}$ of the highest degree
$\deg_1(A_{\si\cup\tau})=|\si|+|\tau|$. Thus the generators
$A_{(1)},A_{(2)},\dots$ are algebraically independent, and the adjoined
graded algebra is naturally isomorphic to the algebra of polynomials in
$A_{(1)},A_{(2)},\dots$ (over the field
$\Bbb C$).

Another filtration on $\ai$ was introduced in~\cite{7}.

\proclaim{Proposition 10.1}
The function
$$
\deg_2(A_\rho) = |\rho| + m_1(\rho)
\tag 10.2
$$
defines a filtration on the algebra $\ai$.
\endproclaim

\demo{Proof} Consider partial permutations
$(d_1,w_1)$, $(d_2,w_2)$ and break the union of their supports
$d_1\cup d_2$ into disjoint parts as follows:
$$
\aligned
d_1\setminus d_2 &= d^{12}_{mf} \cup d^{12}_{ff} \\
d_2\setminus d_1 &= d^{21}_{mf} \cup d^{21}_{ff} \\
d_1 \cap d_2 &= d_{mm} \cup d_{fm} \cup d_{mf} \cup d_{ff}.
\endaligned
$$
The first index equals $f$, if the points of the corresponding domain are
fixed for the permutation $w_1$, and equals $m$, if they are non-fixed. The
second index has a similar sense with respect to the permutation
$w_2$. By definition,
$$
\aligned
\deg_2(w_1) &= |d^{12}_{mf}| + |d_{mm}| + |d_{mf}| +
2\,|d^{12}_{ff}| + 2\,|d_{ff}| + 2\,|d_{fm}| \\
\deg_2(w_2) &= |d^{21}_{fm}| + |d_{mm}| + |d_{fm}| +
2\,|d^{21}_{ff}| + 2\,|d_{ff}| + 2\,|d_{mf}|.
\endaligned
$$

The permutation $w=w_1w_2$ of the set $d=d_1\cup d_2$ has no fixed points in
the domains $d^{12}_{mf}$, $d_{mf}$, $d_{fm}$,
$d^{21}_{fm}$. On the contrary, all points of the domains $d^{12}_{ff}$,
$d_{ff}$, $d^{21}_{ff}$ are fixed for $w$. The domain $d_{mm}$ may contain
both non-fixed and fixed points of $w$.
Hence
$$
\aligned
\deg_2(w_1w_2)
&\le |d^{12}_{mf}| + |d_{mf}| + |d_{fm}| +\\
&+ |d^{21}_{fm}| +
2\,|d^{12}_{ff}| + 2\,|d_{ff}| + 2\,|d^{21}_{ff}| + 2\,|d_{mm}| \\
&\le |d^{12}_{mf}| + 3\,|d_{mf}| + 3\,|d_{fm}| +\\
&+ |d^{21}_{fm}| +
2\,|d^{12}_{ff}| + 4\,|d_{ff}| + 2\,|d^{21}_{ff}| + 2\,|d_{mm}| \\
&= \deg_2(w_1) + \deg_2(w_2),
\endaligned
$$
and the Proposition follows.
\qed\enddemo

Let $\Cal T_2$ be the set of transpositions, i.e. of all permutations from
$\s\infty$ with a unique non-trivial cycle of length~$2$.
The length $\deg(w)$ of a permutation $w$ with respect to a family of
generators $\Cal T_2$ is called the Cayley metric. It is clear that the function
$\deg_3(d,w)=\deg(w)$ defines a filtration on the algebra $\ai$.
The Cayley filtration was studied in~\cite{5}, \cite{9, Chap.~I, \S7,
examples~24,~25}, \cite{6}. In particular, it is known that
$$
\deg_3(A_\rho) = |\rho| - \ell(\rho).
\tag 10.3
$$
Note that according to~\cite{5, Lemma~3.9},
$\deg_3(A_\si)+\deg_3(A_\tau)=\deg_3(A_\rho)$ if and only if the polynomial
$q_{\si,\tau}^\rho(n)$ introduced in Proposition~7.3 is a constant not
depending on $n$.

One may set a problem of general description of filtrations on the algebra
$\ai$. Not having a general answer, we make here only several
observations. First of all, let us give some definitions.

Denote by $\Bbb Y$ the set of all partitions. A function
$\theta:\Bbb Y\to\Bbb Z_+$ is called a {\it filtration} of the algebra
$\ai$, if each triple of partitions $\si,\tau,\rho$ with
$g_{\si,\tau}^\rho>0$ satisfies the inequality
$\theta(\rho)\le\theta(\si)+\theta(\tau)$. We say that a filtration $\theta$
is {\it additive}, if
$$
\theta(\si\cup\tau) = \theta(\si)+\theta(\tau)
\tag 10.4
$$
for all $\si,\tau\in\Bbb Y$. Condition~(10.4) means that
$$
\theta(\rho) = \sum_{k\ge1} \ga_k\;m_k(\rho),
\tag 10.5
$$
where $\ga_k=\theta((k))$ are the degrees of one-cycle permutations.

\example{Example 10.1} All above-mentioned filtrations $\deg_1$,
$\deg_2$, $\deg_3$ are additive. The constants $\gamma_k$ are of the form
$$
\aligned
&\gamma_1=1,\;\gamma_2=2,\;\gamma_3=3,\;\gamma_4=4,\;\ldots
\qquad \text{ for }\quad \deg_1 \\
&\gamma_1=2,\;\gamma_2=2,\;\gamma_3=3,\;\gamma_4=4,\;\ldots
\qquad \text{ for }\quad \deg_2 \\
&\gamma_1=0,\;\gamma_2=1,\;\gamma_3=2,\;\gamma_4=3,\;\ldots
\qquad \text{ for }\quad \deg_3.
\endaligned
$$
\endexample

Let us mention some common properties of the constants $\gamma_k$.

\proclaim{Proposition 10.2}
The following properties of the numbers $\gamma=\{\gamma_k\}_{k=1}^\infty$
are common for all additive filtrations:
$$
\gather
0 \le \gamma_1 \le \gamma_2 \le\gamma_3 \le\gamma_4 \le \ldots
\tag 10.6 \\
\gamma_{i+j+1} \le \gamma_{i+1} + \gamma_{j+1}.
\tag 10.7 \\
k\,\gamma_1 \le 2\,\gamma_k
\tag 10.8 \\
\gamma_{k+1} \le k\,\gamma_2
\tag 10.9 \\
\gamma_{2k+1} \le 2\,\gamma_{k+1}.
\tag 10.10
\endgather
$$
The limit
$$
\lim_{k\to\infty} {\gamma_{k+1}\over k} =
\inf_k {\gamma_{k+1}\over k} =: L(\gamma),
\tag 10.11
$$
always exists, and $\gamma_1\le 2\,L(\gamma)\le2\,\gamma_2$.
\endproclaim

\demo{Proof}
Since $A_{(1)}*A_{(1)}=2A_{(1^2)}+A_{(1)}$, we have
$\gamma_1\le2\gamma_1$, and hence $\gamma_1\ge0$. Multiplying cycles that
intersect by a common pair of neighbour elements, we have
$$
\gathered
(b_1,b_2,\,\ldots,b_i,a_1,a_2)\quad
(a_1,a_2,c_1,c_2,\,\ldots,c_j) = \\ =
(b_1,b_2,\,\ldots,b_i,a_1) (a_2,c_1,c_2,\,\ldots,c_j),
\endgathered
\tag 10.12
$$
which implies that $g_{(i+2),(j+2)}^{(i+1),(j+1)}>0$ and
$\gamma_{i+1}+\gamma_{j+1}\le\gamma_{i+2}+\gamma_{j+2}$. In particular, for
$i=j$ we obtain $\gamma_j\le\gamma_{j+1}$ which proves~(10.6).

To prove~(10.7), note that
$$
(b_1,b_2,\,\ldots,b_i,a)\;
(a,c_1,c_2,\,\ldots,c_j) =
(b_1,b_2,\,\ldots,b_i,a,c_1,c_2,\,\ldots,c_j),
$$
thus $g_{(i+1),(j+1)}^{(i+j+1)}>0$. Statement~(10.11) is a standard
corollary of inequalities~(10.7) (see, for example,~\cite{4, problem~98}).

Formula~(10.8) follows from
$$
(b_1,b_2,\,\ldots,b_k)\;
(b_k,b_{k-1},\,\ldots,b_1) =
(b_1)\,(b_2)\,\ldots\,(b_k),
$$
and~(10.10) is a particular case of~(10.7) for $i=j=k$.

Since
$$
\gathered
(a_1,a_2)\,(a_3,a_4)\,\ldots\,(a_{2k-1},a_{2k})\qquad
(a_2,a_3)\,(a_4,a_5)\,\ldots\,(a_{2k},a_{2k+1}) = \\ =
(a_2,a_4,\,\ldots,a_{2k},a_{2k+1},a_{2k-1},\,\ldots,a_3,a_1); \\
(a_1,a_2)\,(a_3,a_4)\,\ldots\,(a_{2k+1},a_{2k+2})\qquad
(a_2,a_3)\,(a_4,a_5)\,\ldots\,(a_{2k},a_{2k+1}) = \\ =
(a_2,a_4,\,\ldots,a_{2k},a_{2k+2},a_{2k+1},a_{2k-1},\,\ldots,a_3,a_1),
\endgathered
$$
we obtain inequalities~(10.9) for even and odd $k$ respectively. The last
statement of the Proposition follows immediately from the definition of the
limit $L(\gamma)$ and formulae~(10.8),~(10.9).
\qed\enddemo

\remark{Remark}
Conditions (10.6) --- (10.11) do not guarantee that the function defined by
the numbers $\gamma$ via formula~(10.5) is a filtration. For example, the
number of non-trivial cycles
$\theta(\rho)=\ell(\rho)-m_1(\rho)$ corresponds to the constants
$\gamma_1=0$, $\gamma_k=1$ for $k\ge2$. Inequalities (10.6) --- (10.11)
are satisfied, but the function $\theta$ is not a filtration:
$$
(1,2,3,4)\; (1,5,4,6,3)=(1,5)\; (2,3)\; (4,6).
$$
\endremark

\remark{Remark}
Formula~(10.12) may be generalized as follows:
$$
\gathered
(b_1,b_2,\,\ldots,b_i,a_1,a_2,\,\ldots,a_{2k})\quad
(a_1,a_2,\,\ldots,a_{2k},c_1,c_2,\,\ldots,c_j) =\\=
(b_1,b_2,\,\ldots,b_i,a_1,a_3,\,\ldots,a_{2k-1})
(a_2,a_4,\,\ldots,a_{2k},c_1,c_2,\,\ldots,c_j)
\endgathered
$$
for an even number of common elements of multiplied cycles. For an odd
number of common elements we have
$$
\gathered
(b_1,b_2,\ldots,b_i,a_1,a_2,\ldots,a_{2k},a_{2k+1})
(a_1,a_2,\ldots,a_{2k},a_{2k+1},c_1,c_2,\ldots,c_j)\\=
(b_1,b_2,\ldots,b_i,a_1,a_3,\ldots,a_{2k-1},a_{2k+1},
c_1,c_2,\ldots,c_j,a_2,a_4,\ldots,a_{2k}),
\endgathered
$$
thus $\gamma_{i+j+2k+1}\le\gamma_{i+2k+1}+\gamma_{j+2k+1}$.
Note also that
$$
\gathered
(b_1,b_2,\,\ldots,b_i,a_0,a_1,\,\ldots,a_{k-1})\quad
(a_{k-1},\,\ldots,a_1,a_0,c_1,c_2,\,\ldots,c_j) =\\=
(b_1,b_2,\,\ldots,b_i,a_0,c_1,c_2,\,\ldots,c_j)
(a_1)\,(a_2)\,\,\ldots\,(a_{k-1}),
\endgathered
$$
thus
$\gamma_{i+j+1}+(k-1)\gamma_1\le\gamma_{i+k}+\gamma_{j+k}$ for
$k\ge1$.
\endremark

\remark{Remark}
The equality
$$
\gathered
(a_1,\,\ldots,a_i,w,
c_1,\,\ldots,c_k,v,
b_1,\,\ldots,b_j,u) \\
(\alpha_1,\,\ldots,\alpha_I,u,
\beta_1,\,\ldots,\beta_J,v,
\gamma_1,\,\ldots,\gamma_K,w) = \\ =
(a_1,\,\ldots,a_i,w,\alpha_1,\,\ldots,\alpha_I)\,
(b_1,\,\ldots,b_j,u,\beta_1,\,\ldots,\beta_J)\times\\
\times
(c_1,\,\ldots,c_k,v,\gamma_1,\,\ldots,\gamma_K)
\endgathered
$$
implies the inequality
$\gamma_{I+i+1}+\gamma_{J+j+1}+\gamma_{K+k+1}\le
\gamma_{I+J+K+3}+\gamma_{i+j+k+3}$ and, in particular,
$3\gamma_{2n-1}\le2\gamma_{3n}$.
\endremark

\proclaim{Proposition 10.3}
Let $J\subset\Bbb N$. Put
$$
\theta_J(\rho)=|\rho|+\sum_{k\in J}m_k(\rho).
$$
Then the function $\theta_J$ is an additive filtration of the algebra
$\ai$.
\endproclaim

\demo{Proof}
Consider partitions $\si,\tau,\rho$ such that
$g_{\si,\tau}^\rho>0$. Then there exist elements
$(d_1,\om_\si)\in A_\si,\, (d_2,\om_\tau)\in A_\tau$ such that
$(d_1\cup d_2,\om_\si \om_\tau)\in A_\rho$.

Consider the decompositions  of the permutations
$\om_\si,\,\om_\tau,\,\om_\si\om_\tau$ into products of disjoint cycles. If
$1\not\in J$, then denote by
$M$ the set of cycles of the permutation $\om_\si\om_\tau$ that
are contained neither in the decomposition of $\om_\si$, nor in the
decomposition of $\om_\tau$, and
the length of each such cycle belongs to $J$.
If $1\in J$, then we also include in the set $M$ fixed points (informally
speaking, ``cycles of length one'') of the permutation
$\om_\si\om_\tau$ in the set
$d_1\cup d_2$ that  are fixed neither with respect to
$\om_\si$, nor with respect to $\om_\tau$. Then
$$
|M|\ge\sum_{k\in J}(m_k(\rho)-m_k(\si)-m_k(\tau)).
\tag10.13
$$

Let $\alpha$ be a cycle from $M$. Then there exists at least one
element in this cycle that belongs to $d_1\cap d_2$.
In a similar way, if $x$ is a fixed point belonging to $M$, then
$x\in d_1\cap d_2$. Hence
$$
|M|\le |d_1\cap d_2|=|d_1|+|d_2|-|d_1\cup d_2|=|\si|+|\tau|-|\rho|.
\tag10.14
$$

It follows from~(10.13) and~(10.14) that
$$
|\si|+|\tau|+\sum_{k\in J}(m_k(\si)+m_k(\tau))\ge
|\rho|+\sum_{k\in J}m_k(\rho).
$$
\qed\enddemo

\subhead 11. Examples of convolutions of classes in the algebra $\ai$ \endsubhead
We present below the simplest formulae for multiplication of basis elements
$a_\rho$ of the algebra $\ai$. In view of Theorem~9.1, one may regard the same
formulae as examples of multiplication of the functions
$p_\rho^\#$ in the algebra
$\Lambda^*$.

One may also use similar formulae for calculating the convolution of
conjugacy classes in symmetric groups. For example, substituting
$a_\rho=z_\rho\,A_\rho$ we obtain from the corresponding row the formula
$$
A_{(3)} * A_{(3)} =
  2\, A_{(3^2)} +
  5\, A_{(5)} +
  8\, A_{(2^2)} +
  3\, A_{(3 1)} +
      A_{(3)} +
  2\, A_{(1^3)}.
$$
Passing to the homomorphic images $\psi(A_\rho)$ and substituting
$$
\psi(A_\rho) = {t-r+m_1(\rho)\choose m_1(\rho)}\;A_\rho,
$$
we obtain an example of the formula from~\cite{5},
$$
C_{(3)} * C_{(3)} =
  2\, C_{(3^2)} +
  5\, C_{(5)} +
  8\, C_{(2^2)} +
(3t-8)\, C_{(3)} +
{t(t-1)(t-2) \over 3}\, C_\varnothing.
$$
Convolutions of conjugacy classes are obtained from this formula by
substituting different positive integer values of the variable
$t$. For example,
$$
\aligned
C_{(3);3} * C_{(3);3} &=
  2\, C_{\varnothing;3} +
  C_{(3);3} \\
C_{(3);4} * C_{(3);4} &=
  8\, C_{\varnothing;4} +
  4\, C_{(3);4} +
  8\, C_{(2^2);4} \\
C_{(3);5} * C_{(3);5} &=
 20\, C_{\varnothing;5} +
  7\, C_{(3);5} +
  8\, C_{(2^2);5} +
  5\, C_{(5);5} \\
C_{(3);6} * C_{(3);6} &=
 40\, C_{\varnothing;6} +
 10\, C_{(3);6} +
  8\, C_{(2^2);6} +
  5\, C_{(5);6} +
  2\, C_{(3^2);6}\;.
\endaligned
$$

So, examples of multiplication formulae:
$$
\aligned
a_{(2)}\,\,a_{(2)} = &
  a_{(2^2)} +
  4\, a_{(3)} +
  2\, a_{(1^2)}
\endaligned
$$
$$
\aligned
a_{(3)}\,\,a_{(2)} = &
  a_{(3 2)} +
  6\, a_{(4)} +
  6\, a_{(2 1)}
\endaligned
$$
$$
\aligned
a_{(4)}\,\,a_{(2)} = &
  a_{(4 2)} +
  8\, a_{(5)} +
  4\, a_{(2^2)} +
  8\, a_{(3 1)} \\
a_{(2^2)}\,\,a_{(2)} = &
  a_{(2^3)} +
  8\, a_{(3 2)} +
  8\, a_{(4)} +
  4\, a_{(2 1^2)} \\
a_{(3)}\,\,a_{(3)} = &
  a_{(3^2)} +
  9\, a_{(5)} +
  9\, a_{(2^2)} +
  9\, a_{(3 1)} +
  3\, a_{(3)} +
  3\, a_{(1^3)}
\endaligned
$$
$$
\aligned
a_{(5)}\,\,a_{(2)} = &
  a_{(5 2)} +
  10\, a_{(6)} +
  10\, a_{(3 2)} +
  10\, a_{(4 1)} \\
a_{(3 2)}\,\,a_{(2)} = &
  a_{(3 2^2)} +
  6\, a_{(4 2)} +
  4\, a_{(3^2)} +
  12\, a_{(5)} +
  6\, a_{(2^2 1)} +
  2\, a_{(3 1^2)} \\
a_{(4)}\,\,a_{(3)} = &
  a_{(4 3)} +
  12\, a_{(6)} +
  24\, a_{(3 2)} +
  12\, a_{(4 1)} +
  12\, a_{(4)} +
  12\, a_{(2 1^2)} \\
a_{(2^2)}\,\,a_{(3)} = &
  a_{(3 2^2)} +
  12\, a_{(4 2)} +
  24\, a_{(5)} +
  12\, a_{(2^2 1)} +
  24\, a_{(3 1)}
\endaligned
$$
$$
\aligned
a_{(6)}\,\,a_{(2)} &=
  a_{(6 2)} +
  12\, a_{(7)} +
  12\, a_{(4 2)} +
  6\, a_{(3^2)} +
  12\, a_{(5 1)} \\
a_{(42)}\,\,a_{(2)} &=
  a_{(4 2^2)} +
  8\, a_{(5 2)} +
  4\, a_{(4 3)} +
  16\, a_{(6)} +
  4\, a_{(2^3)} +\\
&+  8\, a_{(3 2 1)} +
  2\, a_{(4 1^2)} \\
a_{(3^2)}\,\,a_{(2)} &=
  a_{(3^2 2)} +
  12\, a_{(4 3)} +
  18\, a_{(6)} +
  12\, a_{(3 2 1)} \\
a_{(2^3)}\,\,a_{(2)}  &=
  a_{(2^4)} +
  12\, a_{(3 2^2)} +
  24\, a_{(4 2)} +
  6\, a_{(2^2 1^2)}
\endaligned
$$
$$
\aligned
a_{(5)}\,\,a_{(3)} &=
  a_{(5 3)} +
  15\, a_{(7)} +
  30\, a_{(4 2)} +
  15\, a_{(3^2)} +\\
  &+15\, a_{(5 1)} +
  30\, a_{(5)} +
   15\, a_{(2^2 1)} +
  15\, a_{(3 1^2)} \\
a_{(3 2)}\,\,a_{(3)} &=
  a_{(3^2 2)} +
  9\, a_{(5 2)} +
  6\, a_{(4 3)} +
  36\, a_{(6)} +\\
  &+9\, a_{(2^3)} +
  15\, a_{(3 2 1)}
  + 21\, a_{(3 2)} +
  36\, a_{(4 1)} +
  3\, a_{(2 1^3)} \\
a_{(4)}\,\,a_{(4)} &=
  a_{(4^2)} +
  16\, a_{(7)} +
  32\, a_{(4 2)} +
  24\, a_{(3^2)} +
  16\, a_{(5 1)}
  +48\, a_{(5)} +\\
  &+32\, a_{(2^2 1)} +
  + 4\, a_{(2^2)} +
  16\, a_{(3 1^2)} +
  16\, a_{(3 1)} +
  4\, a_{(1^4)}
\endaligned
$$
$$
\aligned
a_{(2^2)}\,\,a_{(2^2)} &=
  a_{(2^4)} +
  16\, a_{(3 2^2)} +
  32\, a_{(4 2)} +
  32\, a_{(3^2)} +\\
  &+64\, a_{(5)} +
  8\, a_{(2^2 1^2)} +
  32\, a_{(3 1^2)} +
 + 16\, a_{(2^2)} +
  8\, a_{(1^4)}
\endaligned
$$

\subhead 12. Irreducible representations of the semigroup $\Cal P_n$ and
characters of the algebra $\bi$ \endsubhead
Fix a subset $x\subset\Bbb P_n$ of size $|x|=k$ and let
$\lambda$ be a Young diagram with $k$ boxes. Then the formula
$$
\pi_{x,\lambda}(d,w) = \cases
\pi_\lambda(\varphi_x(d,w)), & \text{ if } d \subset x; \\
0 & \text{ otherwise}
\endcases
\tag 12.1
$$
defines an irreducible representation of the semigroup
$\p n$. Here $0$ is the zero matrix of order
$\dim\pi_\lambda$.

\proclaim{Proposition 12.1}
The representations $\pi_{x,\lambda}$ (where $x\subset\Bbb P_n$, $|x|=k$
and $\lambda\in\Bbb Y_k$ are a subset and a Young diagram of common size
$0\le k\le n$) are irreducible, pairwise non-equivalent and form a
complete
list of irreducible representations of the semigroup algebra
$\b n$.
\endproclaim

\demo{Proof}
The matrices $\pi_{x,\lambda}(d,w)$ of the representation $\pi_{x,\lambda}$
are non-zero exactly for elements
$(d,w)$ such that
$d\subset x$. Thus the representations indexed by different sets
$x$ are non-equivalent. It is obvious that the representations
$\pi_{x,\lambda}$, $\pi_{x,\mu}$ for $\lambda\ne\mu$ are non-equivalent too.
Completeness of the list follows from Corollary~3.2.
\qed\enddemo

Let us describe the branching rule for irreducible representations of the
semigroup
$\p n$ when restricting on $\Cal P_{n-1}\subset\p n$.

We write
$\mu\nearrow\lambda$,
if a diagram $\la$
is obtained from a diagram
$\mu$ by adding one box. Let
$\Gamma_n$ be the set of pairs $(x,\lambda)$, where
$x\subset\Bbb P_n$ is a subset, $\lambda$ is a Young diagram, and
$x$ and $\lambda$ are assumed to have the same size,
$|x|=|\lambda|$. The set $\Gamma_n$ indexes irreducible representations
$\pi_{x,\lambda}$ of the semigroup $\p n$.

\proclaim{Proposition 12.2}
Assume that a set  $x\subset\Bbb P_n$ does not contain
$n$. Then the restriction of the representation $\pi_{x,\lambda}$ on
the subsemigroup $\Cal P_{n-1}$ remains irreducible (and is indexed by the
same pair $(x,\lambda)$ regarded as an element of
$\Gamma_{n-1}$). And if $n\in x$, then the restriction of the irreducible
representation $\pi_{x,\lambda}$ on $\Cal P_{n-1}$ is of the form
$$
\res \pi_{x,\lambda} = \bigoplus_{\mu:\mu\nearrow\lambda}
\pi_{y,\mu},
\tag 12.2
$$
where $y=x\setminus\{n\}$.
\endproclaim

\demo{Proof} Follows immediately from the construction of representations
and a well known branching rule for irreducible representations of symmetric
groups.
\qed\enddemo

Thus the set of vertices of the branching graph of irreducible
representations of the semigroups $\p n$ is
$\Gamma=\bigcup_{n=0}^\infty\Gamma_n$. Let $(y,\mu)$,
$(x,\lambda)$ be vertices of neighbour levels $\Gamma_{n-1}$ and
$\Gamma_n$. They are joined by an edge (of multiplicity~$1$) if and only
if $\mu=\lambda$ or $\mu\nearrow\lambda$ (i.e. if the diagram
$\lambda$ coincides with $\mu$ or is obtained from $\mu$
by adding one box). In the first case $x=y$, in the second case
$x=y\cup\{n\}$.

Applying the ergodic method (see~\cite{2}) to the branching graph
$\Gamma$, one can easily obtain a description of non-negative indecomposable
harmonic functions on $\Gamma$ or, equivalently, of characters of the
algebra $\bi$. Not considering this question in details, we point
out only the parameterization of characters.

Let $X$ be an arbitrary subset in $\Bbb
P_\infty=\{1,2,\ldots\}$. Denote by $\Gamma_X$ the set of vertices
$(x,\lambda)\in\Gamma$ with $x\subset X$.
By Proposition~12.2, $\Gamma_X$ is a coideal in the branching graph
$\Gamma$, i.e.
$$
(y,\mu) \nearrow (x,\lambda) \in \Gamma_X \quad \Longrightarrow
\quad (y,\mu) \in \Gamma_X.
$$

Each harmonic function on a coideal can be canonically continued (by zero)
to a harmonic function on the whole branching graph. Note that if the set
$X$ is infinite, then the graph $\Gamma_X$
differs from the Young graph $\Bbb Y$ only by trivial doubling of some
levels. As for the Young graph, its indecomposable harmonic functions
$\varphi^{(\alpha;\beta;X)}$ are indexed by the points of the Thoma simplex
$\Delta$. By definition,
$\Delta$ consists of pairs of non-increasing non-negative sequences
$\alpha=(\alpha_1,\alpha_2,\ldots)$,
$\beta=(\beta_1,\beta_2,\ldots)$ with
$$
\sum_{n=1}^\infty \alpha_n +
\sum_{n=1}^\infty \beta_n \le 1.
$$
The functions $\varphi^{(\alpha;\beta;X)}$ are of the form
$$
\varphi^{(\alpha;\beta;X)}(y,\mu) = \cases
s_\mu(\alpha;\beta),& \text{if } y\subset X \\
0,& \text{otherwise},
\endcases
$$
where $s_\mu(\alpha;\beta)$ are the extended Schur functions,
see~\cite{1}.

If $X$ is a finite set, $|X|=k$, then indecomposable harmonic functions
$\varphi^{(\lambda;X)}$ on the coideal
$\Gamma_X$ are indexed by Young diagrams $\lambda\in\Bbb Y_k$
with $k$ boxes and are of the form
$$
\varphi^{(\lambda;X)}(y,\mu) = \cases
{\dim(\mu,\lambda) \over \dim\lambda}, & \text{if }
y\subset X\; \text{and }\; \mu\subset\lambda \\
0, & \text{otherwise},
\endcases
$$
where $\dim(\mu,\lambda)$ is the number of standard skew Young tableaux of
shape $\lambda/\mu$ and $\dim\lambda$ is the number of all standard
tableaux of shape $\lambda$.

Translated by Natalia Tsilevich.

\enddocument